\documentclass[11pt,reqno]{amsart}
\usepackage[a4paper,left=29mm,right=29mm,top=27mm,bottom=28mm]{geometry}
\usepackage[T1]{fontenc}
\usepackage{lmodern}
\usepackage{amsmath,amssymb,amsthm}
\usepackage{enumitem}
\usepackage{needspace}
\usepackage{microtype}
\usepackage[colorlinks=true,linkcolor=blue,citecolor=blue,urlcolor=blue]{hyperref}
\hypersetup{
	pdftitle={Optimal quantitative stability near the Coulomb fission threshold},
	pdfauthor={},
	pdfsubject={Sharp first-order degeneration and almost optimal configurations}
}

\numberwithin{equation}{section}
\newtheorem{theorem}{Theorem}[section]
\newtheorem{proposition}[theorem]{Proposition}
\newtheorem{lemma}[theorem]{Lemma}
\newtheorem{corollary}[theorem]{Corollary}
\theoremstyle{remark}
\newtheorem{remark}[theorem]{Remark}
\setlist[enumerate]{leftmargin=2.2em,label=\textup{(\roman*)},itemsep=0.35em,topsep=0.4em}
\title[Quantitative stability at the fission threshold]
{Optimal quantitative stability near the Coulomb fission threshold}
\author{Zijuan Gao, Chengxiang Zhang}
\address[]{Laboratory of Mathematics and Complex Systems (Ministry of Education), 
	School of Mathematical Sciences, 
	Beijing Normal University, Beijing 100875, P. R. China}
\email{gaozijuan3012@163.com}	
\address[]{Laboratory of Mathematics and Complex Systems (Ministry of Education), 
	School of Mathematical Sciences, 
	Beijing Normal University, Beijing 100875, P. R. China}
\email{zcx@bnu.edu.cn}	
\date{\today}

\begin{document}
	\begin{abstract}
		We study the optimal quadratic stability constant for the three-dimensional Coulomb liquid-drop energy at fixed volume, measured by the squared Fraenkel asymmetry. Let $M$ be the critical mass at which the energy of one ball equals the limiting energy of two infinitely separated balls of equal volume. We prove that, as $m$ approaches $M$ from below, the optimal constant satisfies $\Lambda(m)=(36\pi)^{1/3}(2^{1/3}-1)(1-m/M)+o(1-m/M)$. We further show that every almost optimal sequence with masses $m_n\to M$ from below is asymptotic in measure to a union of two equal-volume balls whose centers diverge. We also identify the corresponding perimeter and Coulomb-energy limits and give a necessary and sufficient criterion for almost optimality after rescaling to mass $M$. The proof combines a strict binding inequality for mass partitions, a local spectral gap at the critical ball, and a BV truncation argument for critical minimizing sequences. All results hold for arbitrary finite-perimeter competitors, without connectedness or boundary regularity assumptions.
	\end{abstract}
		\maketitle

	\small{
		}
	\smallskip
	
	\section{Introduction}
	
	Let $E\subset\mathbb R^3$ be a measurable set of finite volume and finite perimeter. We consider
\begin{equation}\label{eq:energy}
	\mathcal E_1(E)=P(E)+D(E),
	\qquad
	D(E)=\frac12\iint_{E\times E}\frac{dx\,dy}{|x-y|}.
\end{equation}
The perimeter $P(E)$ is taken in the sense of
De Giorgi, namely
\[
P(E)
=
\sup\left\{
\int_E \operatorname{div}F(x)\,dx
\;\middle|\;
F\in C_c^1(\mathbb R^3;\mathbb R^3),\
|F(x)|\leq 1\ \text{for all }x\in\mathbb R^3
\right\},
\]
	which is simply the surface area of $E$ when the boundary
	is smooth.
	
	This functional is the three-dimensional Coulomb liquid-drop energy introduced in Gamow's model of atomic nuclei. The perimeter term represents surface tension, whereas the Coulomb term describes the repulsive interaction between protons. We refer to Choksi, Muratov, and Topaloglu \cite{CMT} for a survey of the model. A related problem with a background potential was studied by Alama, Bronsard, Choksi, and Topaloglu \cite{ABCT}.
	
	For each $m>0$ we consider the minimal energy at fixed mass,
	\[
	\inf_{\substack{|E|=m,\;P(E)<\infty}}\mathcal E_1(E).
	\]
	For $0<m<M$, where $M$ is the critical mass defined in
	\eqref{eq:critical-mass}, the infimum is attained by a ball. This paper concerns
	the quantitative stability of that minimizer and its sharp
	degeneration as $m\to M^-$.
	
	The two terms in \eqref{eq:energy} have competing effects. At fixed volume, balls minimize the perimeter and maximize the Coulomb energy. For $|E|=1$, the scaling relation is
	\[
	\mathcal E_1(m^{1/3}E)=m^{2/3}\bigl(P(E)+mD(E)\bigr).
	\]
	Thus the perimeter dominates at small masses, whereas the Coulomb interaction becomes more significant at large masses.
	
	Write $B_m=B_{R_m}(0)$, where $R_m=(3m/(4\pi))^{1/3}$, and set $a=(36\pi)^{1/3}$. The ball energy is
	\[
	f(m)=\mathcal E_1(B_m)=a m^{2/3}+\frac a5m^{5/3}.
	\]
	The critical fission mass $M$ is determined by $f(M)=2f(M/2)$ and is given by
	\begin{equation}\label{eq:critical-mass}
		M=5\frac{2^{1/3}-1}{1-2^{-2/3}}
		=5\frac{2-2^{2/3}}{2^{2/3}-1}.
	\end{equation}
	
	The global minimization problem for the liquid-drop energy has been studied extensively. Knüpfer and Muratov \cite{KM13,KM14} established existence results for small masses and nonexistence results for sufficiently large masses, together with spherical rigidity in several regimes. Further results on spherical symmetry and stability were obtained by Julin \cite{Julin} and by Figalli, Fusco, Maggi, Millot, and Morini \cite{FFMMM}. In the three-dimensional Coulomb case, Chodosh and Ruohoniemi \cite{CR} proved spherical uniqueness for $m\leq1$, while Frank, Killip, and Nam \cite{FKN} obtained further nonexistence results for large masses. For general Riesz interactions, Frank and Nam \cite{FN} proved existence up to the corresponding critical mass and established a conditional uniqueness result.
	
	The closest analytic precedents for our argument are the compactness theory of Frank and Lieb \cite{FL} and the second-variation criterion of Bonacini and Cristoferi \cite{BC}. Frank and Lieb developed a $BV$ compactness framework for the liquid-drop model, including a missing-mass argument. Bonacini and Cristoferi showed that positivity of the second variation yields a quantitative local minimality inequality in the $L^1$ topology for a class of nonlocal isoperimetric problems. Related work concerns second-variation stability for attractive-repulsive energies \cite{BCT}, large-volume nonlocal isoperimetric problems under confinement \cite{GO}, the dilute limit of nuclear matter \cite{EFK}, and the sharp-interface limit connecting the TFDW and liquid-drop models \cite{ASB}. These works provide relevant analytic tools but address different interaction kernels, background fields, or asymptotic regimes. Here we consider the pure Coulomb energy, without confinement or an external background potential. We use the criterion in \cite{BC} for local stability and a $BV$ truncation argument to analyze critical minimizing sequences.
	
	Chodosh and Gianocca \cite{CG} established the global-minimizer classification for the three-dimensional Coulomb liquid-drop model: balls are the unique minimizers, up to translation, for $0<m\leq M$, while no minimizer exists for $m>M$. We use its endpoint assertion at $m=M$, recorded in Proposition~\ref{prop:critical-input}. This classification concerns actual minimizers. It does not, however, determine the sharp degeneration of the optimal stability constant or the structure of all almost optimal sequences as $m\to M^-$.
	
	This leaves two quantitative questions. The first is the exact first-order behavior of the optimal quadratic stability constant as $m\to M^-$; the second is the possible limiting configurations of sequences that attain this optimal constant asymptotically. The latter is linked to the loss of compactness at the critical mass: two balls of volume $M/2$, placed at a distance tending to infinity, form a minimizing sequence for the problem at mass $M$. Thus the endpoint is governed not only by the ball minimizer but also by configurations that split into distant components.
	
	For $|E|=m>0$, the Fraenkel asymmetry is
	\begin{equation}\label{eq:asymmetry}
		\mathcal A(E)=\frac1m\inf_{x\in\mathbb R^3}|E\triangle(x+B_m)|.
	\end{equation}
	In particular, $\mathcal A(E)=1$ for two equal-volume balls with sufficiently separated centers. For $0<m<M$, define
	\begin{equation}\label{eq:optimal-constant}
		\Lambda(m)=
		\inf_{\substack{|E|=m,\;P(E)<\infty\\ \mathcal A(E)>0}}
		\frac{\mathcal E_1(E)-f(m)}{m^{2/3}\mathcal A(E)^2}.
	\end{equation}
	
	Our first result answers the first question: it determines the exact
	first-order behavior of $\Lambda(m)$ at the fission threshold and
	supplies two-sided bounds that hold at every subcritical mass.
	
	\begin{theorem}\label{thm:main}
		Let $\Lambda(m)$ be defined by \eqref{eq:optimal-constant}. Then
		\begin{equation}\label{eq:sharp-limit}
			\lim_{m\to M^-}\frac{\Lambda(m)}{1-m/M}=a(2^{1/3}-1),
		\end{equation}
		where $a=(36\pi)^{1/3}$. Moreover, there exists a constant $c_{\mathrm{iso}}>0$ such that for all $0<m<M$,
		\begin{equation}\label{eq:coarse-bounds}
			c_{\mathrm{iso}}(1-m/M)\leq\Lambda(m)
			\leq a(2^{1/3}-1)(1-m/M).
		\end{equation}
	\end{theorem}
	
	The limit in \eqref{eq:sharp-limit} shows that the global stability constant vanishes linearly at the fission threshold. This degeneration is not a loss of local stability: the ball remains strictly stable under sufficiently small volume-preserving perturbations at mass $M$. Rather, it is caused by competitors that split into distant components. The leading coefficient has a direct geometric interpretation as the normalized perimeter cost of splitting a ball of mass $M$ into two equal balls:
	\[
	a(2^{1/3}-1)=\frac{2P(B_{M/2})-P(B_M)}{M^{2/3}}.
	\]
	
	The upper bound follows from the exact two-ball comparison in Lemma~\ref{lem:two-balls}. For the lower bound, the scaling identity separates the perimeter deficit from the critical energy excess. After rescaling competitors to mass $M$, the local spectral gap excludes convergence to a single ball. A $BV$ truncation then produces two separated profiles; the strict binding inequality forces their masses to converge to $M/2$, and the quantitative isoperimetric inequality makes both profiles asymptotically spherical. The perimeter term in the scaling identity yields the sharp lower bound.
	
	Our second main result gives the complete structure of almost optimal sequences.
	\begin{theorem}\label{thm:almost-optimal}
		Let $m_n\to M$ with $0<m_n<M$, and set $\delta_n=1-m_n/M$. Let $E_n\subset\mathbb R^3$ satisfy $|E_n|=m_n$, $P(E_n)<\infty$, $\mathcal A(E_n)>0$, and
		\begin{equation}\label{eq:almost-optimal}
			\frac{\mathcal E_1(E_n)-f(m_n)}
			{m_n^{2/3}\mathcal A(E_n)^2}
			=\Lambda(m_n)+o(\delta_n).
		\end{equation}
		Then there exist $x_n^+,x_n^-\in\mathbb R^3$ such that
		\begin{equation}\label{eq:two-ball-convergence}
			|x_n^+-x_n^-|\longrightarrow\infty,
			\qquad
			\frac{|E_n\triangle((x_n^++B_{m_n/2})\cup(x_n^-+B_{m_n/2}))|}{m_n}
			\longrightarrow0.
		\end{equation}
		Moreover,
		\begin{equation}\label{eq:component-limits}
			\mathcal A(E_n)\to1,\qquad
			\frac{P(E_n)}{m_n^{2/3}}\to a2^{1/3},\qquad
			\frac{D(E_n)}{m_n^{5/3}}\to\frac a5\,2^{-2/3}.
		\end{equation}
		The dilated sets $F_n=(M/m_n)^{1/3}E_n$ satisfy
		\begin{equation}\label{eq:critical-small-excess}
			\mathcal E_1(F_n)-f(M)=o(\delta_n).
		\end{equation}
	\end{theorem}
	
	Theorem~\ref{thm:almost-optimal} shows that the two-ball splitting is necessary for every almost optimal sequence, not merely for the explicit recovery family used to prove the upper bound. In particular, the conclusion holds for arbitrary finite-perimeter competitors and requires neither connectedness nor boundary regularity.
	
	The paper is organized as follows. Section~\ref{sec:scaling} establishes the scaling identity and the exact two-ball comparison. Section~\ref{sec:binding-local} proves the strict binding inequality and the local stability estimate at the critical ball. Section~\ref{sec:compactness} establishes the $BV$ profile decomposition and classifies critical minimizing sequences. Section~\ref{sec:main-proofs} proves the sharp asymptotic formula, characterizes almost optimal sequences, and gives an equivalent criterion for almost optimality.
	
	\Needspace{8\baselineskip}
	\section{Scaling and two-ball comparison}\label{sec:scaling}
	We first record the potential and scaling formulas. We then compute the interaction between two separated balls and obtain the upper bound in \eqref{eq:coarse-bounds}.
	The Newtonian potential of a ball is classical; we record the formulas
	and refer to \cite[Lemma~2]{CR} for the computation.
	For $R>0$ set $V_R(x)=\int_{B_R(0)}|x-y|^{-1}\,dy$ and write $r=|x|$. Then
	\begin{equation}\label{eq:ball-potential}
		V_R(x)=
		\begin{cases}
			2\pi(R^2-r^2/3),&0\leq r\leq R,\\
			4\pi R^3/(3r),&r\geq R,
		\end{cases}
	\end{equation}
	and $D(B_R)=\frac{16\pi^2}{15}R^5$. Since $P(B_R)=4\pi R^2$ while
	$m=|B_R|=\frac43\pi R^3$, we have $D(B_m)=\frac m5P(B_m)$, and hence
	\begin{equation}\label{eq:ball-energy}
		P(B_m)=a m^{2/3},\qquad
		D(B_m)=\frac a5m^{5/3},\qquad
		f(m)=a m^{2/3}+\frac a5m^{5/3}.
	\end{equation}
	
	\subsection{Continuity and asymmetry}
	
	\begin{lemma}\label{lem:continuity}
		Let $E,F\subset\mathbb R^3$ satisfy $|E|,|F|\leq K<\infty$. Then
		\begin{equation}\label{eq:coulomb-continuity}
			|D(E)-D(F)|\leq(2\pi+K)|E\triangle F|.
		\end{equation}
		If $|E|=|F|=m>0$, then
		\begin{equation}\label{eq:asymmetry-continuity}
			|\mathcal A(E)-\mathcal A(F)|
			\leq\frac{|E\triangle F|}{m}.
		\end{equation}
		For sets of positive volume, $\mathcal A(tE)=\mathcal A(E)$ for $t>0$, and
		$0\leq\mathcal A(E)\leq2$.
	\end{lemma}
	
	\begin{proof}
		For every $x\in\mathbb R^3$,
		\[
		\int_E\frac{dy}{|x-y|}
		\leq\int_{B_1(x)}\frac{dy}{|x-y|}
		+|E|
		\leq2\pi+K.
		\]
		Thus every finite-volume set has finite Coulomb energy.
		Set $V_E(x)=\int_E|x-y|^{-1}\,dy$. By symmetry of the kernel, we obtain
		\[
		D(E)-D(F)
		=\frac12\int_{\mathbb R^3}
		(\mathbf1_E-\mathbf1_F)(V_E+V_F)\,dx.
		\]
		
		It follows that
		\[
		|D(E)-D(F)|
		\leq\frac{\|V_E\|_\infty+\|V_F\|_\infty}{2}
		|E\triangle F|
		\leq(2\pi+K)|E\triangle F|.
		\]
		
		For any $x\in\mathbb R^3$, the triangle inequality for symmetric
		difference implies
		\[
		|E\triangle(x+B_m)|
		\leq|E\triangle F|+|F\triangle(x+B_m)|.
		\]
		Taking the infimum and interchanging $E$ and $F$ yields
		\eqref{eq:asymmetry-continuity}.
		
		For $|E|=m$, definition~\eqref{eq:asymmetry} gives
		\[
		\frac{|tE\triangle(x+B_{t^3m})|}{t^3m}
		=\frac{|E\triangle(x/t+B_m)|}{m}.
		\]
		Taking the infimum over $x$ proves scale invariance.
		
		Finally, $|E\triangle(x+B_m)|\leq2m$ yields
		$\mathcal A(E)\leq2$.
	\end{proof}
	
	\subsection{The critical global input and the scaling identity}
	
	We use the following established result of Chodosh and Gianocca,
	with the normalization in \eqref{eq:energy}.
	
	\begin{proposition}\label{prop:critical-input}
		Let $F\subset\mathbb R^3$ be a finite-perimeter set with $|F|=M$. Then
		\[
		\mathcal E_1(F)\geq f(M).
		\]
		Equality holds if and only if $F$ is a translate of $B_M$.
	\end{proposition}
	
	This is the endpoint case of \cite[Theorem~1.1]{CG}.
	Only this case of that classification theorem is needed below.
	
	\begin{lemma}\label{lem:scaling}
		Let $0<m\leq M$ and let $E\subset\mathbb R^3$ satisfy
		$|E|=m$ and $P(E)<\infty$. Set $t=(m/M)^{1/3}$ and $F=t^{-1}E$.
		Then $|F|=M$ and
		\begin{equation}\label{eq:scaling-identity}
			\frac{\mathcal E_1(E)-f(m)}{m^{2/3}}
			=
			\frac{(1-t^3)(P(F)-P(B_M))
				+t^3(\mathcal E_1(F)-f(M))}{M^{2/3}}.
		\end{equation}
		In particular,
		\begin{equation}\label{eq:perimeter-lower}
			\mathcal E_1(E)-f(m)
			\geq(1-m/M)\bigl(P(E)-P(B_m)\bigr).
		\end{equation}
	\end{lemma}
	
	\begin{proof}
		
		For $X\in C_c^1(\mathbb R^3;\mathbb R^3)$,
		$\operatorname{div}_u(X(tu))=t(\operatorname{div}X)(tu)$.
		The change $Y(u)=X(tu)$ is bijective on the vector fields of
		supremum norm at most one. Hence
		\begin{equation}\label{eq:perimeter-scaling}
			\begin{aligned}
				P(tF)
				=\,\sup_{\|X\|_\infty\leq1}
				\int_{tF}\operatorname{div}X\,dx
				=\,t^2\sup_{\|X\|_\infty\leq1}
				\int_F\operatorname{div}_u(X(tu))\,du
				=t^2P(F),
			\end{aligned}
		\end{equation}
		and
		\begin{equation}\label{eq:coulomb-scaling}
			D(tF)=\frac{t^6}{2}
			\iint_{F\times F}\frac{du\,dv}{t|u-v|}=t^5D(F).
		\end{equation}
		
		By \eqref{eq:perimeter-scaling} and
		\eqref{eq:coulomb-scaling}, together with $B_m=tB_M$,
		we obtain
		\[
		\begin{aligned}
			\mathcal E_1(E)-f(m)
			=&\,t^2(P(F)-P(B_M))+t^5(D(F)-D(B_M))\\
			=&\,t^2\bigl((1-t^3)(P(F)-P(B_M))
			+t^3(\mathcal E_1(F)-f(M))\bigr).
		\end{aligned}
		\]
		
		Dividing by $m^{2/3}=t^2M^{2/3}$ gives
		\eqref{eq:scaling-identity}.
		By Proposition~\ref{prop:critical-input},
		$\mathcal E_1(F)-f(M)\geq0$. Hence
		\eqref{eq:scaling-identity} yields
		\[
		\mathcal E_1(E)-f(m)
		\geq t^2(1-t^3)\bigl(P(F)-P(B_M)\bigr).
		\]
		Since $t^3=m/M$ and
		$P(E)-P(B_m)=t^2(P(F)-P(B_M))$,
		we obtain \eqref{eq:perimeter-lower}.
		
	\end{proof}
	
	\begin{corollary}\label{cor:subcritical}
		There exists a dimensional constant $c_{\mathrm{iso}}>0$ such that
		\begin{equation}\label{eq:subcritical-stability}
			\mathcal E_1(E)-f(m)
			\geq c_{\mathrm{iso}}(1-m/M)m^{2/3}\mathcal A(E)^2
			\quad\text{for }0<m<M,\quad |E|=m.
		\end{equation}
		For each such mass, balls are the unique minimizers up to translation.
	\end{corollary}
	
	\begin{proof}
		By the quantitative isoperimetric inequality
		\cite{FMP,CL}, we get
		\[
		P(E)-P(B_m)\geq
		c_{\mathrm{iso}}m^{2/3}\mathcal A(E)^2.
		\]
		Combining this with \eqref{eq:perimeter-lower} proves
		\eqref{eq:subcritical-stability}.
		If $\mathcal E_1(E)=f(m)$, then
		\eqref{eq:perimeter-lower} and $m<M$ imply $P(E)=P(B_m)$.
		Hence, by the equality case of the isoperimetric inequality,
		$E$ is a translate of $B_m$.
	\end{proof}
	
	\subsection{Exact interaction between two balls}
	\begin{lemma}\label{lem:two-balls}
		Let $m>0$, $R=R_m$, $r=R_{m/2}$, and $L>2(R+r)$.
		Set $E_L=B_r(0)\cup B_r(Le_1)$, where
		$L$ denotes the distance between the centers of the two balls.
		Then
		\begin{equation}\label{eq:two-ball-exact}
			\mathcal A(E_L)=1,\qquad
			\mathcal E_1(E_L)-f(m)
			=a(2^{1/3}-1)m^{2/3}(1-m/M)+\frac{m^2}{4L}.
		\end{equation}
	\end{lemma}

	\begin{proof}
		Since $L>2(R+r)$, the two components are disjoint.
		If a translate of $B_m$ intersected both components,
		the triangle inequality would imply $L\leq2(R+r)$,
		a contradiction. Thus every translate of $B_m$
		intersects $E_L$ in volume at most $m/2$.
		The ball $B_R(0)$ contains $B_r(0)$ and is disjoint
		from $B_r(Le_1)$. Hence
		\[
		\mathcal A(E_L)
		=2-\frac2m\sup_{x\in\mathbb R^3}
		|E_L\cap(x+B_m)|
		=2-\frac2m\frac m2=1.
		\]
		
		By \eqref{eq:ball-potential}, for $x\in B_r(0)$, there holds
		\[
		\int_{B_r(Le_1)}\frac{dy}{|x-y|}
		=\frac{m}{2|x-Le_1|}.
		\]
		Since $x\mapsto |x-Le_1|^{-1}$ is harmonic
		on a neighborhood of $\overline{B_r(0)}$,
		the mean-value property gives
		\begin{equation}\label{eq:two-ball-cross-interaction}
			\begin{aligned}
				\iint_{B_r(0)\times B_r(Le_1)}
				\frac{dx\,dy}{|x-y|}
				&=\frac m2\int_{B_r(0)}
				\frac{dx}{|x-Le_1|}
				=\frac{m^2}{4L}.
			\end{aligned}
		\end{equation}
		The two cross terms in $D(E_L)$ are equal, so the
		factor $1/2$ leaves one copy of this integral.
		Consequently,
		\[
		\mathcal E_1(E_L)=2f(m/2)+\frac{m^2}{4L}.
		\]
		
		By \eqref{eq:ball-energy} and \eqref{eq:critical-mass},
		\[
		\begin{aligned}
			\mathcal E_1(E_L)-f(m)
			=&\,2f(m/2)-f(m)+\frac{m^2}{4L}
			=\,a(2^{1/3}-1)m^{2/3}(1-m/M)
			+\frac{m^2}{4L}.
		\end{aligned}
		\]
		This proves \eqref{eq:two-ball-exact}.
	\end{proof}

	Corollary~\ref{cor:subcritical} yields the lower bound
	in \eqref{eq:coarse-bounds}, while Lemma~\ref{lem:two-balls}
	yields the upper bound as $L\to\infty$.
	Thus $\Lambda(m)>0$ for $m<M$, but the sharp limit
	\eqref{eq:sharp-limit} requires a further argument.

	\section{Binding and local stability at the critical mass}\label{sec:binding-local}
	
	\subsection{Strict binding and its equality cases}\label{sec:binding}

	The following scalar inequality identifies the mass partitions at the critical mass that can attain the minimum energy. It is a special case of the binding comparisons in \cite{FN}; we determine all equality cases explicitly.

	\begin{lemma}\label{lem:binding}
		Let $x,y>0$ satisfy $x+y\leq M$. Then
		\begin{equation}\label{eq:scalar-binding}
			f(x)+f(y)\geq f(x+y).
		\end{equation}
		If $x+y<M$, the inequality is strict. If $x+y=M$, equality holds
		if and only if $x=y=M/2$.
	\end{lemma}
	\begin{proof}
		Set $m=x+y$, $\theta=x/m$, and
		\[
		h(\theta)=
		\frac{\theta^{2/3}+(1-\theta)^{2/3}-1}
		{1-\theta^{5/3}-(1-\theta)^{5/3}},
		\qquad 0<\theta<1.
		\]
		The denominator is positive since $s^{5/3}<s$
		for $0<s<1$. We show that $h$ has a unique minimum
		at $\theta=1/2$.
		
		Let $u=\theta^{1/3}$, $v=(1-\theta)^{1/3}$,
		and $w=u+v$. Then
		\[
		u^3+v^3=1,\qquad uv=\frac{w^3-1}{3w},
		\qquad 1<w\leq2^{2/3},
		\]
		where equality in the upper bound holds only when
		$u=v$, by the strict concavity of $s^{1/3}$.
		Using
		$u^5+v^5=(u^2+v^2)(u^3+v^3)-u^2v^2(u+v)$,
		we obtain
		\[
		u^2+v^2-1
		=\frac{(w-1)^2(w+2)}{3w},
		\]
		and
		\[
		\begin{aligned}
			1-u^5-v^5
			=&\,1-(u^2+v^2)+(uv)^2w
			=\,\frac{9w-3(w^3+2)+(w^3-1)^2}{9w}\\
			=&\,\frac{(w-1)^3(w^3+3w^2+6w+5)}{9w}.
		\end{aligned}
		\]
		Consequently,
		\[
		h(\theta)
		=\frac{3(w+2)}
		{(w-1)(w^3+3w^2+6w+5)}.
		\]
		Differentiation gives
		\[
		\frac{d}{dw}h
		=-\frac{9(w^2+w+1)(w^2+3w+1)}
		{(w-1)^2(w^3+3w^2+6w+5)^2}<0.
		\]
		Hence $h$ is strictly decreasing in $w$. Since
		$w\leq2^{2/3}$, with equality only at $\theta=1/2$,
		we obtain
		\begin{equation}\label{eq:binding-ratio}
			h(\theta)\geq h(1/2)
			=\frac{2^{1/3}-1}{1-2^{-2/3}}
			=\frac M5,
		\end{equation}
		with equality only at $\theta=1/2$, where the last
		identity follows from \eqref{eq:critical-mass}.
		
		By \eqref{eq:ball-energy}, there holds
		\[
		f(x)+f(y)-f(m)
		=a m^{2/3}
		\bigl(1-\theta^{5/3}-(1-\theta)^{5/3}\bigr)
		\left(h(\theta)-\frac m5\right).
		\]
		For $m<M$, \eqref{eq:binding-ratio} gives
		\[
		h(\theta)-\frac m5\geq\frac{M-m}{5}>0.
		\]
		For $m=M$, equality holds precisely when
		$\theta=1/2$, by \eqref{eq:binding-ratio}.
		This proves the claimed inequality and its equality case.
	\end{proof}

	\begin{corollary}\label{cor:partitions}
		Let $J\geq1$, $m_j>0$, and $\sum_{j=1}^Jm_j=M$. Then
		\begin{equation}\label{eq:partition-energy}
			\sum_{j=1}^J f(m_j)\geq f(M).
		\end{equation}
		Equality holds if and only if $J=1$ and $m_1=M$, or $J=2$ and
		$m_1=m_2=M/2$.
	\end{corollary}
	
	\begin{proof}
		The case $J=1$ is immediate, while the case $J=2$
		follows from Lemma~\ref{lem:binding}.
		Suppose $J\geq3$. Since all masses are positive
		and sum to $M$, we have $m_1+m_2<M$.
		Thus the first merge strictly decreases the sum
		of the ball energies by \eqref{eq:scalar-binding}.
		Each subsequent merge does not increase this sum,
		again by \eqref{eq:scalar-binding}, since every
		intermediate total is at most $M$.
		The final sum is $f(M)$. Hence
		\eqref{eq:partition-energy} is strict for $J\geq3$.
	\end{proof}
	
	\subsection{The local spectral gap at the critical ball}\label{sec:local}

	The critical ball remains strictly stable under small volume-preserving perturbations, even though the global stability constant in Theorem~\ref{thm:main} vanishes as $m\to M^-$. We compute the constrained spherical spectrum and apply the second-variation criterion of \cite{AFM,BC} in the normalization \eqref{eq:energy}.

	\begin{lemma}\label{lem:local-stability}
		There exist $\eta\in(0,1)$ and $\kappa>0$ such that every set
		$F\subset\mathbb R^3$ with $|F|=M$, $P(F)<\infty$, and
		$\mathcal A(F)<\eta$ satisfies
		\begin{equation}\label{eq:critical-local-stability}
			\mathcal E_1(F)-f(M)
			\geq\kappa M^{2/3}\mathcal A(F)^2.
		\end{equation}
	\end{lemma}

	\begin{proof}
		Set $R=R_M$. By \eqref{eq:ball-potential},
		the restriction of $V_R$ to $\partial B_R(0)$ is constant,
		while its normal derivative satisfies
		\[
		\partial_\nu V_R(x)
		=-\frac{4\pi R}{3}
		=-\frac{M}{R^2},
		\qquad x\in\partial B_R(0).
		\]
		Since the mean curvature of $\partial B_R(0)$ is constant,
		so is the Euler--Lagrange quantity $H+V_R$.
		Thus $B_R(0)$ is a regular volume-constrained critical set.
		
		Let $\varphi\in H^1(\mathbb S^2)$ denote the scalar normal
		velocity, with value $\varphi(\omega)$ at $R\omega$.
		The linearized volume constraint and translation
		orthogonality are
		\begin{equation}\label{eq:normal-constraints}
			\int_{\mathbb S^2}\varphi\,d\omega=0,
			\qquad
			\int_{\mathbb S^2}\varphi(\omega)\omega_i\,d\omega=0,
			\quad i=1,2,3.
		\end{equation}
		By \cite[Theorem~3.3 and equation~(3.2)]{BC},
		with their coefficient $\gamma=1/2$, the constrained
		second variation is
		\begin{equation}\label{eq:spherical-hessian}
			\begin{aligned}
				\mathcal Q_R(\varphi)
				=&\,\underbrace{
					\int_{\mathbb S^2}
					\bigl(|\nabla_{\mathbb S^2}\varphi|^2-2\varphi^2\bigr)
					\,d\omega}_{I_1}
				+\underbrace{
					R^3\iint_{\mathbb S^2\times\mathbb S^2}
					\frac{\varphi(\omega)\varphi(\eta)}{|\omega-\eta|}
					\,d\omega\,d\eta}_{I_2}
				+\underbrace{
					\left(-M\int_{\mathbb S^2}\varphi^2\,d\omega\right)
				}_{I_3}.
			\end{aligned}
		\end{equation}
		
		Define $\widetilde\varphi\in H^1(\partial B_R(0))$ by
		\[
		\widetilde\varphi(R\omega)=\varphi(\omega),
		\qquad \omega\in\mathbb S^2.
		\]
		Since the outward normal satisfies $\nu(R\omega)=\omega$,
		the constraints in \eqref{eq:normal-constraints} identify
		$\widetilde\varphi$ with an element of
		$T^\perp(\partial B_R(0))$. This is the subspace denoted by
		$T^\perp$ in \cite[Definition~3.8]{BC}.
		We compute $I_1$, $I_2$, and $I_3$ separately:
		
		\[
		\begin{aligned}
			I_1
			=&\,\int_{\partial B_R(0)}
			\left(|\nabla_\tau\widetilde\varphi|^2
			-\frac2{R^2}\widetilde\varphi^2\right)\,dS,\quad
			I_2
			=\,\iint_{\partial B_R(0)\times\partial B_R(0)}
			\frac{\widetilde\varphi(x)\widetilde\varphi(y)}
			{|x-y|}\,dS_x\,dS_y,\\[2mm]
			I_3
			=&\,\int_{\partial B_R(0)}
			\partial_\nu V_R\,\widetilde\varphi^2\,dS.
		\end{aligned}
		\]

		Let $\{Y_{\ell j}\}_{j=1}^{2\ell+1}$ be a real orthonormal
		basis of spherical harmonics of degree $\ell$.
		By the Funk--Hecke formula, the eigenvalue of the Coulomb
		operator on the space of spherical harmonics of degree
		$\ell$ is
		\[
		2\pi\int_{-1}^1
		\frac{P_\ell(t)}{\sqrt{2-2t}}\,dt,
		\]
		where $P_\ell$ is the Legendre polynomial of degree
		$\ell$, normalized by $P_\ell(1)=1$. By Rodrigues' formula and direct differentiation,
		
		\[
		P_\ell(t)=\frac1{2^\ell\ell!}
		\frac{d^\ell}{dt^\ell}(t^2-1)^\ell,\qquad
		\frac{d^\ell}{dt^\ell}(1-t)^{-1/2}
		=\frac{\Gamma(\ell+1/2)}{\Gamma(1/2)}
		(1-t)^{-\ell-1/2}.
		\]

		For $\ell\geq1$ and $0\leq j\leq\ell-1$, the boundary
		products satisfy
		\begin{equation}\label{eq:legendre-boundary}
			\begin{aligned}
				\frac{d^j}{dt^j}(t^2-1)^\ell
				\frac{d^{\ell-1-j}}{dt^{\ell-1-j}}(1-t)^{-1/2}
				=&\,O((1-t)^{1/2})\quad(t\uparrow1),\\
				\frac{d^j}{dt^j}(t^2-1)^\ell
				\frac{d^{\ell-1-j}}{dt^{\ell-1-j}}(1-t)^{-1/2}
				=&\,O((1+t)^{\ell-j})\quad(t\downarrow-1).
			\end{aligned}
		\end{equation}
		The estimates in \eqref{eq:legendre-boundary} show that
		all boundary terms vanish. Integration by parts $\ell$ times
		and the beta integral yield
		\[
		\begin{aligned}
			\int_{-1}^1\frac{P_\ell(t)}{\sqrt{2-2t}}\,dt
			=&\,\frac{\Gamma(\ell+1/2)}
			{\sqrt{2\pi}\,2^\ell\ell!}
			\int_{-1}^1(1-t)^{-1/2}(1+t)^\ell\,dt\\
			=&\,\frac{\Gamma(\ell+1/2)}
			{\sqrt{2\pi}\,2^\ell\ell!}
			2^{\ell+1/2}
			\frac{\Gamma(\ell+1)\Gamma(1/2)}
			{\Gamma(\ell+3/2)}
			=\frac{2}{2\ell+1}.
		\end{aligned}
		\]
		For $\ell=0$, the beta integral gives the same value.
		Thus, for every $\ell\geq0$, the Funk--Hecke eigenvalue is
		\begin{equation}\label{eq:coulomb-eigenvalue}
			2\pi\int_{-1}^1
			\frac{P_\ell(t)}{\sqrt{2-2t}}\,dt
			=\frac{4\pi}{2\ell+1}.
		\end{equation}
		
		The constraints in \eqref{eq:normal-constraints} eliminate
		the spherical harmonics of degrees zero and one. Thus
		\[
		\varphi=\sum_{\ell=2}^\infty\sum_{j=1}^{2\ell+1}
		c_{\ell j}Y_{\ell j}.
		\]
		Using \eqref{eq:coulomb-eigenvalue} and $4\pi R^3=3M$,
		substitution into \eqref{eq:spherical-hessian} gives
		\begin{equation}\label{eq:hessian-spectrum}
			\mathcal Q_R(\varphi)
			=
			\sum_{\ell=2}^\infty\sum_{j=1}^{2\ell+1}
			(\ell-1)\left(\ell+2-\frac{2M}{2\ell+1}\right)c_{\ell j}^2.
		\end{equation}
		Since finite spherical-harmonic sums are dense in
		$H^1(\mathbb S^2)$ and the Coulomb operator is bounded
		on $L^2(\mathbb S^2)$ with operator norm $4\pi$,
		the identity extends to all $\varphi\in H^1(\mathbb S^2)$
		satisfying \eqref{eq:normal-constraints}.
		
		For $\ell\geq2$, $(\ell+2)(2\ell+1)\geq20$. Hence,
		for $M<10$,
		\[
		(\ell-1)\left(\ell+2-\frac{2M}{2\ell+1}\right)
		\geq(1-M/10)\bigl(\ell(\ell+1)-2\bigr)
		\geq\frac47(1-M/10)\bigl(1+\ell(\ell+1)\bigr).
		\]
		The last inequality follows from $\ell(\ell+1)\geq6$.
		With
		$\|\varphi\|_{H^1(\mathbb S^2)}^2
		=\int_{\mathbb S^2}(\varphi^2+|\nabla_{\mathbb S^2}\varphi|^2)\,d\omega$,
		we obtain
		\begin{equation}\label{eq:hessian-coercivity}
			\mathcal Q_R(\varphi)
			\geq\frac47(1-M/10)
			\|\varphi\|_{H^1(\mathbb S^2)}^2
			\qquad(M<10).
		\end{equation}
		By \eqref{eq:critical-mass},
		$M=5\,2^{2/3}/(1+2^{1/3})<5\,2^{1/3}<10$.
		Then \eqref{eq:hessian-coercivity} implies that
		$\mathcal Q_{R_M}$ is coercive on the subspace of
		$H^1(\mathbb S^2)$ satisfying \eqref{eq:normal-constraints}.
		
		We apply the local quantitative minimality theorem
		\cite[Theorem~2.8]{BC} to $B_M$. The translation distance
		in this theorem is
		\[
		\inf_{x\in\mathbb R^3}|F\triangle(x+B_M)|
		=M\mathcal A(F).
		\]
		The Coulomb kernel $K(x)=|x|^{-1}$ satisfies the
		Riesz-kernel hypothesis, since $0<1<N-1=2$ for $N=3$.
		Moreover, $B_M$ is smooth and volume-constrained critical.
		By the identification above,
		\eqref{eq:hessian-coercivity} verifies the
		positive-second-variation condition of Bonacini and Cristoferi
		\cite[Definition~3.8]{BC}.
		Consequently, there exist $\eta>0$ and $C>0$ such that
		\[
		\mathcal E_1(F)-f(M)\geq C M^2\mathcal A(F)^2
		=\kappa M^{2/3}\mathcal A(F)^2,
		\qquad \kappa=C M^{4/3}>0.
		\]
		The estimate holds for every finite-perimeter set $F$
		with $|F|=M$ and $0<\mathcal A(F)<\eta$.
		No graph representation of $\partial F$ is required.
		For $\mathcal A(F)=0$, the estimate follows from
		Proposition~\ref{prop:critical-input}.
		We decrease $\eta$ if necessary to ensure
		$\eta\in(0,1)$.
	\end{proof}
	
	\begin{corollary}\label{cor:uniform-local}
		Let $\eta$ and $\kappa$ be as in Lemma~\ref{lem:local-stability}.
		If $M/2\leq m\leq M$, $|E|=m$, and $0<\mathcal A(E)<\eta$, then
		\begin{equation}\label{eq:uniform-local-quotient}
			\frac{\mathcal E_1(E)-f(m)}
			{m^{2/3}\mathcal A(E)^2}\geq\frac\kappa2.
		\end{equation}
	\end{corollary}
	
	\begin{proof}
		Use $F=(M/m)^{1/3}E$ in \eqref{eq:scaling-identity}.
		The perimeter deficit is nonnegative and
		$\mathcal A(F)=\mathcal A(E)$. Consequently,
		\[
		\frac{\mathcal E_1(E)-f(m)}
		{m^{2/3}\mathcal A(E)^2}
		\geq\frac mM
		\frac{\mathcal E_1(F)-f(M)}
		{M^{2/3}\mathcal A(F)^2}
		\geq\frac mM\,\kappa\geq\frac\kappa2.
		\]
	\end{proof}
	
	\section{Critical minimizing sequences}\label{sec:compactness}

	We use the concentration function from the concentration-compactness principle of \cite{Lions} and the $BV$ compactness method of Frank and Lieb \cite{FL}. The required truncation estimates follow from the coarea formula and the slicing theory for sets of finite perimeter; we use the results in Chapters~12--16 and 18 of \cite{Maggi}.
	For open $U\subset\mathbb R^3$, write
	$P(E;U)=|D\mathbf1_E|(U)$.
	\begin{lemma}\label{lem:profile-input}
		
		Let $F_n\subset\mathbb R^3$ satisfy $|F_n|=M$ and
		$\sup_nP(F_n)<\infty$. After passing to a subsequence, one of the
		following alternatives holds.
		\begin{enumerate}
			\item There exist $z_n\in\mathbb R^3$ and a finite-perimeter set $F$
			with $|F|=M$ such that $|(F_n-z_n)\triangle F|\to0$.
			\item There exist $\alpha\in(0,M)$ and disjoint finite-perimeter
			sets $A_n,B_n\subset F_n$ such that
			\begin{equation}\label{eq:profile-mass}
				|A_n|\to\alpha,\qquad |B_n|\to M-\alpha,\qquad
				\operatorname{dist}(A_n,B_n)\to\infty,
			\end{equation}
			and
			\begin{equation}\label{eq:profile-approximation}
				|F_n\setminus(A_n\cup B_n)|\leq\frac4n,\qquad
				P(A_n)+P(B_n)\leq P(F_n)+\frac8{n^2}.
			\end{equation}
			In particular,
			\begin{equation}\label{eq:truncation-energy}
				\mathcal E_1(A_n)+\mathcal E_1(B_n)
				\leq\mathcal E_1(F_n)+\frac8{n^2}.
			\end{equation}
		\end{enumerate}
		The bounds involving $n$ refer to the extracted, relabeled sequence.
		
	\end{lemma}
	
	\begin{proof}
		
		Set
		\[
		q_n(r)=\sup_{x\in\mathbb R^3}|F_n\cap B_r(x)|.
		\]
		Since $0\leq q_n(k)\leq M$, a diagonal argument yields,
		after passing to a subsequence, 
		\[
		q_n(k)\longrightarrow q(k),
		\qquad\text{for every }k\in\mathbb N.
		\]
		The limit $q$ is nondecreasing and bounded by $M$. Thus
		\[
		\alpha:=\lim_{k\to\infty}q(k)\in[0,M]
		\]
		is well defined.
		
		Assume first that $\alpha=0$. For $z\in\mathbb Z^3$, let
		$Q_z=z+(0,1)^3$ and set $v_{n,z}=|F_n\cap Q_z|$.
		Since each $Q_z$ is contained in a ball of radius $2$,
		$v_{n,z}\leq q_n(2)\to0$. Thus $v_{n,z}\leq1/2$ for all
		$z$ and sufficiently large $n$. The relative isoperimetric
		inequality yields
		\[
		v_{n,z}^{2/3}\leq C P(F_n;Q_z).
		\]
		Since the cubes have pairwise disjoint interiors and their
		boundaries are Lebesgue-null, we obtain
		\[
		\begin{aligned}
			M=\sum_zv_{n,z}
			\leq&\,\bigl(\sup_zv_{n,z}\bigr)^{1/3}
			\sum_zv_{n,z}^{2/3}\\
			\leq&\,C q_n(2)^{1/3}\sum_zP(F_n;Q_z)
			\leq C q_n(2)^{1/3}P(F_n)\longrightarrow0,
		\end{aligned}
		\]
		a contradiction. Hence $\alpha>0$.
		
		Assume $\alpha=M$. Fix an integer $r_0$ such that
		$q(r_0)>3M/4$. Since $q_n(r_0)\to q(r_0)$, for all
		sufficiently large $n$, we can choose $z_n$ such that
		\[
		|F_n\cap B_{r_0}(z_n)|>M/2.
		\]
		By $BV$ compactness, after passing to a subsequence,
		\[
		\mathbf1_{F_n-z_n}\longrightarrow\mathbf1_F
		\qquad\text{in }L^1_{\mathrm{loc}}(\mathbb R^3).
		\]
		We fix this subsequence before choosing $\varepsilon$.
		Fix $0<\varepsilon<M/4$ and choose an integer
		$r_\varepsilon$ such that
		$q(r_\varepsilon)>M-\varepsilon/2$.
		For all sufficiently large $n$, choose $y_n$ such that
		\[
		|F_n\cap B_{r_\varepsilon}(y_n)|>M-2\varepsilon.
		\]
		Since
		\[
		|F_n\cap B_{r_0}(z_n)|
		+|F_n\cap B_{r_\varepsilon}(y_n)|
		>|F_n|=M,
		\]
		the sets $F_n\cap B_{r_0}(z_n)$ and
		$F_n\cap B_{r_\varepsilon}(y_n)$ have an intersection
		of positive measure. Hence
		$B_{r_0}(z_n)\cap B_{r_\varepsilon}(y_n)\neq\emptyset$. Consequently,
		\begin{equation}\label{eq:concentration-tightness}
			|y_n-z_n|\leq r_0+r_\varepsilon,\qquad
			|F_n\setminus B_{r_0+2r_\varepsilon}(z_n)|
			<2\varepsilon.
		\end{equation}
		Set $R=r_0+2r_\varepsilon$. The tail estimate
		\eqref{eq:concentration-tightness} and local
		convergence give $|F\cap B_R|\geq M-2\varepsilon$.
		Moreover, local convergence implies $|F|\leq M$.
		Since $\varepsilon$ is arbitrary, $|F|=M$.
		
		It follows that
		\[
		\begin{aligned}
			|(F_n-z_n)\triangle F|
			\leq&\,|((F_n-z_n)\triangle F)\cap B_R|
			+|(F_n-z_n)\setminus B_R|+|F\setminus B_R|\\
			\leq&\,|((F_n-z_n)\triangle F)\cap B_R|
			+4\varepsilon.
		\end{aligned}
		\]
		For each fixed $\varepsilon$, the first term tends to zero
		as $n\to\infty$. Sending $\varepsilon\to0$ then proves
		\[
		\mathbf1_{F_n-z_n}\longrightarrow\mathbf1_F
		\qquad\text{in }L^1(\mathbb R^3).
		\]
		This proves (i).
		
		It remains to consider $0<\alpha<M$. Choose integers $r_j\geq j$
		such that $q(r_j)>\alpha-j^{-1}$. Choose strictly increasing indices $\nu_j\uparrow\infty$
		such that
		\[
		|q_{\nu_j}(r_j)-q(r_j)|
		+|q_{\nu_j}(4r_j)-q(4r_j)|<j^{-1}.
		\]
		For simplicity, we relabel $(F_{\nu_j},r_j)$ as $(F_n,r_n)$. Then
		\[
		r_n\geq n,\qquad
		q_n(r_n)>\alpha-\frac2n,\qquad
		q_n(4r_n)<\alpha+\frac1n.
		\]
		Choose $x_n$ with
		$|F_n\cap B_{r_n}(x_n)|>q_n(r_n)-n^{-1}$. Then
		\[
		\alpha-\frac3n
		<|F_n\cap B_{r_n}(x_n)|
		\leq|F_n\cap B_{4r_n}(x_n)|
		<\alpha+\frac1n,
		\]
		so that
		\begin{equation}\label{eq:concentration-annulus}
			|F_n\cap(B_{4r_n}(x_n)\setminus B_{r_n}(x_n))|
			\leq\frac4n.
		\end{equation}
		Let $F_n^{(1)}$ denote the points of density one.
		By \eqref{eq:concentration-annulus} and coarea on $(r_n,2r_n)$ and $(3r_n,4r_n)$, there exist radii
		$\rho_n\in(r_n,2r_n)$ and $\sigma_n\in(3r_n,4r_n)$ such that
		\begin{equation}\label{eq:cutting-slice-bound}
			\mathcal H^2(F_n^{(1)}\cap\partial B_{\rho_n}(x_n))
			+\mathcal H^2(F_n^{(1)}\cap\partial B_{\sigma_n}(x_n))
			\leq\frac8{nr_n}.
		\end{equation}
		The radii may be chosen among those for which the $BV$ slicing
		identities hold. Set
		\[
		A_n=F_n\cap B_{\rho_n}(x_n),\qquad
		B_n=F_n\setminus\overline{B_{\sigma_n}(x_n)}.
		\]
		For these representatives,
		\[
		\begin{gathered}
			\operatorname{dist}(A_n,B_n)\geq\sigma_n-\rho_n\geq r_n\to\infty,\quad
			\alpha-\frac3n\leq|A_n|\leq\alpha+\frac1n,\\
			M-\alpha-\frac1n\leq|B_n|\leq M-\alpha+\frac3n,\quad
			|F_n\setminus(A_n\cup B_n)|\leq\frac4n.
		\end{gathered}
		\]
		By \eqref{eq:cutting-slice-bound},
		\begin{equation}\label{eq:cutting-perimeter}
			\begin{aligned}
				P(A_n)+P(B_n)
				=&\,P(F_n;B_{\rho_n}(x_n))
				+P(F_n;\mathbb R^3\setminus\overline{B_{\sigma_n}(x_n)})\\
				&+\mathcal H^2(F_n^{(1)}\cap\partial B_{\rho_n}(x_n))
				+\mathcal H^2(F_n^{(1)}\cap\partial B_{\sigma_n}(x_n))\\
				\leq&\,P(F_n)+\frac8{nr_n}
				\leq P(F_n)+\frac8{n^2}.
			\end{aligned}
		\end{equation}
		This proves \eqref{eq:profile-mass}--\eqref{eq:profile-approximation}.
		Since the Coulomb kernel is nonnegative and $A_n,B_n$ are disjoint
		subsets of $F_n$,
		\[
		D(A_n)+D(B_n)\leq D(F_n).
		\]
		Adding this to the perimeter estimate
		\eqref{eq:cutting-perimeter} proves
		\eqref{eq:truncation-energy}.
		
	\end{proof}
	
	The next proposition treats critical minimizing sequences,
	including unbounded sets with arbitrarily many components.
	The truncations above identify their limiting profiles.
	Continuity of the Coulomb term and convergence of the total energy
	then yield the separate Coulomb and perimeter limits, in particular
	the perimeter convergence used in \eqref{eq:sequence-liminf}.
	
	\Needspace{9\baselineskip}
	\begin{proposition}\label{prop:critical-dichotomy}
		Let $F_n\subset\mathbb R^3$ satisfy $|F_n|=M$ and
		$\mathcal E_1(F_n)\to f(M)$. Then, after passing to a subsequence,
		exactly one of the following alternatives holds.
		\begin{enumerate}
			\item There exist $z_n\in\mathbb R^3$ such that
			$|F_n\triangle(z_n+B_M)|\to0$. In this case,
			$\mathcal A(F_n)\to0$, $P(F_n)\to P(B_M)$, and
			$D(F_n)\to D(B_M)$.
			\item There exist $z_n^+,z_n^-\in\mathbb R^3$ with
			$|z_n^+-z_n^-|\to\infty$ such that
			\[
			\left|F_n\triangle
			\bigl((z_n^++B_{M/2})\cup(z_n^-+B_{M/2})\bigr)\right|\to0.
			\]
			In this case,
			\begin{equation}\label{eq:critical-pair-limits}
				\mathcal A(F_n)\to1,\qquad
				P(F_n)\to2P(B_{M/2}),\qquad
				D(F_n)\to2D(B_{M/2}).
			\end{equation}
		\end{enumerate}
	\end{proposition}
	
	\begin{proof}
		
		Since $P(F_n)\leq\mathcal E_1(F_n)$, Lemma~\ref{lem:profile-input}
		applies. In the compactness alternative, there exist
		$z_n\in\mathbb R^3$ and a set $F$ with $|F|=M$ such that
		\[
		|(F_n-z_n)\triangle F|\to0.
		\]
		By Lemma~\ref{lem:continuity} and lower semicontinuity of perimeter, there holds
		\[
		\mathcal E_1(F)
		\leq\liminf_{n\to\infty}\mathcal E_1(F_n)
		=f(M).
		\]
		Proposition~\ref{prop:critical-input} gives the reverse inequality,
		with equality only for translates of $B_M$. Thus $F=z+B_M$
		for some $z\in\mathbb R^3$. Replacing $z_n$ by $z_n+z$, we obtain
		\[
		|(F_n-z_n)\triangle B_M|\to0.
		\]
		By Lemma~\ref{lem:continuity},
		$D(F_n)\to D(B_M)$. Moreover,
		$|(F_n-z_n)\triangle B_M|\to0$ implies
		$\mathcal A(F_n)\to0$. Hence
		\[
		P(F_n)=\mathcal E_1(F_n)-D(F_n)\to P(B_M).
		\]
		
		In the second alternative of Lemma~\ref{lem:profile-input}, set
		$a_n=|A_n|$ and $b_n=|B_n|$. For all sufficiently large $n$,
		$a_n,b_n\in(0,M)$. By Corollary~\ref{cor:subcritical} and
		\eqref{eq:truncation-energy},
		\[
		f(a_n)+f(b_n)
		\leq\mathcal E_1(A_n)+\mathcal E_1(B_n)
		\leq\mathcal E_1(F_n)+\frac8{n^2}.
		\]
		Passing to the limit and applying
		Corollary~\ref{cor:partitions} with $J=2$, we obtain
		\[
		f(M)\leq f(\alpha)+f(M-\alpha)\leq f(M).
		\]
		The equality case in Corollary~\ref{cor:partitions} yields $\alpha=M/2$. Moreover,
		\begin{equation}\label{eq:split-excess}
			\begin{aligned}
				0\leq
				\bigl(\mathcal E_1(A_n)-f(a_n)\bigr)
				+\bigl(\mathcal E_1(B_n)-f(b_n)\bigr)
				\leq\
				\mathcal E_1(F_n)+\frac8{n^2}
				-f(a_n)-f(b_n)
				\longrightarrow0.
			\end{aligned}
		\end{equation}
		For all sufficiently large $n$, $a_n,b_n\in[M/4,3M/4]$.
		Combining \eqref{eq:split-excess} with
		\eqref{eq:subcritical-stability} yields
		\begin{equation}\label{eq:component-asymmetry}
			\mathcal A(A_n)^2+\mathcal A(B_n)^2
			\leq\frac{4^{5/3}}{c_{\mathrm{iso}}M^{2/3}}
			\bigl(\mathcal E_1(A_n)-f(a_n)
			+\mathcal E_1(B_n)-f(b_n)\bigr)\longrightarrow0.
		\end{equation}
		By \eqref{eq:component-asymmetry} and
		\eqref{eq:asymmetry}, choose $z_n^\pm\in\mathbb R^3$ such that
		\begin{equation}\label{eq:component-ball-approximation}
			\begin{aligned}
				|A_n\triangle(z_n^++B_{a_n})|
				\leq a_n\mathcal A(A_n)+n^{-1},\quad
				|B_n\triangle(z_n^-+B_{b_n})|
				\leq b_n\mathcal A(B_n)+n^{-1}.
			\end{aligned}
		\end{equation}
		By \eqref{eq:component-ball-approximation}, both
		intersections with the corresponding balls have positive
		measure for large $n$. Hence
		\begin{equation}\label{eq:component-separation}
			|z_n^+-z_n^-|
			\geq\operatorname{dist}(A_n,B_n)-2R_{3M/4}
			\longrightarrow\infty.
		\end{equation}
		Since concentric balls are nested,
		$|B_{a_n}\triangle B_{M/2}|=|a_n-M/2|$, and similarly
		for $b_n$. Set
		\[
		G_n=(z_n^++B_{M/2})\cup(z_n^-+B_{M/2}).
		\]
		Then, by \eqref{eq:profile-approximation}
		and \eqref{eq:component-ball-approximation},
		\begin{equation}\label{eq:two-ball-profile}
			\begin{aligned}
				|F_n\triangle G_n|
				\leq\,&\frac6n+a_n\mathcal A(A_n)+b_n\mathcal A(B_n)
				+|a_n-M/2|+|b_n-M/2|
				\longrightarrow0.
			\end{aligned}
		\end{equation}
		By \eqref{eq:component-separation}, the two balls are
		disjoint and satisfy $|z_n^+-z_n^-|>2(R_M+R_{M/2})$ for large $n$. Lemma~\ref{lem:two-balls} gives
		$\mathcal A(G_n)=1$. Applying
		\eqref{eq:two-ball-cross-interaction} with $m=M$ and
		$L=|z_n^+-z_n^-|$ yields
		\[
		D(G_n)=2D(B_{M/2})+\frac{M^2}{4|z_n^+-z_n^-|}.
		\]
		By \eqref{eq:two-ball-profile},
		\eqref{eq:coulomb-continuity}, and
		\eqref{eq:asymmetry-continuity},
		\[
		\mathcal A(F_n)\to1,\qquad
		D(F_n)\to2D(B_{M/2}).
		\]
		Since $\mathcal E_1(F_n)\to f(M)$, we also have
		\[
		\begin{aligned}
			P(F_n)
			\to f(M)-2D(B_{M/2})
			=2f(M/2)-2D(B_{M/2})
			=2P(B_{M/2}).
		\end{aligned}
		\]
		This proves \eqref{eq:critical-pair-limits}.
		The two alternatives are mutually exclusive, since their
		Fraenkel asymmetries converge to $0$ and $1$, respectively.
		
	\end{proof}
	
	\begin{remark}
		
		The finite-profile energy in \cite[Lemma~6]{ABCT} has
		Coulomb coefficient $1$. Define
		\[
		\widetilde{\mathcal E}(G)
		=P(G)+\iint_{G\times G}|x-y|^{-1}\,dx\,dy.
		\]
		The dilation $F\mapsto 2^{-1/3}F$ gives
		\[
		\widetilde{\mathcal E}(2^{-1/3}F)
		=2^{-2/3}\mathcal E_1(F),
		\qquad |2^{-1/3}F|=|F|/2.
		\]
		Consequently,
		\[
		\inf_{|G|=M/2}\widetilde{\mathcal E}(G)
		=2^{-2/3}f(M).
		\]
		This identifies the normalization relating the
		decomposition in \cite{ABCT} to the present functional.
		Our proof instead uses the explicit two-piece truncation
		\eqref{eq:profile-approximation}.
		
	\end{remark}
	
	\section{The sharp coefficient and almost optimal sequences}
	\label{sec:main-proofs}
	
	The next proposition establishes the compactness needed for the sharp lower bound under a hypothesis weaker than almost optimality.
	
	\begin{proposition}\label{prop:bounded-quotients}
		Let $m_n\to M$ with $0<m_n<M$, and set $\delta_n=1-m_n/M$.
		Let $E_n\subset\mathbb R^3$ satisfy $|E_n|=m_n$, $P(E_n)<\infty$,
		$\mathcal A(E_n)>0$, and, for some $C<\infty$ independent of $n$,
		\begin{equation}\label{eq:bounded-quotient}
			\frac{\mathcal E_1(E_n)-f(m_n)}
			{m_n^{2/3}\mathcal A(E_n)^2}\leq C\delta_n.
		\end{equation}
		Then \eqref{eq:two-ball-convergence} and
		\eqref{eq:component-limits} hold. Moreover,
		\begin{equation}\label{eq:sequence-liminf}
			\liminf_{n\to\infty}
			\frac{\mathcal E_1(E_n)-f(m_n)}
			{\delta_n m_n^{2/3}\mathcal A(E_n)^2}
			\geq a(2^{1/3}-1).
		\end{equation}
	\end{proposition}
	
	\begin{proof}
		Set $t_n=(m_n/M)^{1/3}$ and $F_n=t_n^{-1}E_n$.
		By \eqref{eq:scaling-identity},
		\begin{equation}\label{eq:quotient-decomposition}
			\begin{aligned}
				\frac{\mathcal E_1(E_n)-f(m_n)}
				{\delta_n m_n^{2/3}\mathcal A(E_n)^2}
				=&\,\frac{P(F_n)-P(B_M)}
				{M^{2/3}\mathcal A(F_n)^2}
				+\frac{t_n^3}{\delta_n}
				\frac{\mathcal E_1(F_n)-f(M)}
				{M^{2/3}\mathcal A(F_n)^2}.
			\end{aligned}
		\end{equation}
		Both terms on the right are nonnegative. Since
		$\mathcal A(F_n)\leq2$, \eqref{eq:bounded-quotient} gives
		\[
		0\leq\mathcal E_1(F_n)-f(M)
		\leq\frac{4CM^{2/3}}{t_n^3}\,\delta_n.
		\]
		Since $t_n^3=m_n/M\to1$, we have $t_n^3\geq1/2$
		for all sufficiently large $n$. Hence
		\begin{equation}\label{eq:rescaled-energy-bound}
			0\leq\mathcal E_1(F_n)-f(M)
			\leq8CM^{2/3}\delta_n\longrightarrow0.
		\end{equation}
		By \eqref{eq:rescaled-energy-bound}, $(F_n)$ is a minimizing sequence at mass $M$.
		
		We next exclude the single-ball alternative without imposing
		any regularity on $E_n$. If $0<\mathcal A(E_n)<\eta$ along a subsequence, then
		\eqref{eq:uniform-local-quotient} and
		\eqref{eq:bounded-quotient} yield
		\[
		\frac{\kappa}{2}\leq C\delta_n,
		\]
		a contradiction to $\delta_n\to0$.
		Consequently,
		\begin{equation}\label{eq:away-from-ball}
			\mathcal A(F_n)=\mathcal A(E_n)\geq\eta,
			\quad\text{for all sufficiently large }n.
		\end{equation}
		
		By \eqref{eq:away-from-ball}, every subsequence has a
		further subsequence satisfying the second alternative of
		Proposition~\ref{prop:critical-dichotomy}. Hence, along each such subsequence,
		\[
		\mathcal A(F_n)\to1,\quad
		P(F_n)\to2P(B_{M/2}),\quad
		D(F_n)\to2D(B_{M/2}).
		\]
		Thus the convergences hold for the full sequence. Since $m_n=t_n^3M$, we obtain
		\[
		\frac{P(E_n)}{m_n^{2/3}}=\frac{P(F_n)}{M^{2/3}}
		\longrightarrow\frac{2a(M/2)^{2/3}}{M^{2/3}}
		=a2^{1/3},
		\]
		and
		\[
		\frac{D(E_n)}{m_n^{5/3}}=\frac{D(F_n)}{M^{5/3}}
		\longrightarrow\frac{(2a/5)(M/2)^{5/3}}{M^{5/3}}
		=\frac a5\,2^{-2/3}.
		\]
		The scale invariance of $\mathcal A$ completes
		\eqref{eq:component-limits}.
		
		Fix $j\geq1$. The preceding double-ball extraction implies
		that, for all sufficiently large $n$, there exist
		$z_{n,j}^\pm\in\mathbb R^3$ such that
		\begin{equation}\label{eq:diagonal-two-ball-approximation}
			|z_{n,j}^+-z_{n,j}^-|\geq j,\qquad
			\left|F_n\triangle
			\bigl((z_{n,j}^++B_{M/2})\cup
			(z_{n,j}^-+B_{M/2})\bigr)\right|\leq\frac1j.
		\end{equation}
		Otherwise, a subsequence would contradict the preceding
		extraction result. Choose strictly increasing integers $N_j$
		such that \eqref{eq:diagonal-two-ball-approximation} holds for all $n\geq N_j$.
		For $N_j\leq n<N_{j+1}$, select the corresponding centers
		and denote them by $z_n^\pm$. Then
		\[
		|z_n^+-z_n^-|\to\infty,\qquad
		|F_n\triangle G_n|\to0,
		\]
		where
		\[
		G_n=(z_n^++B_{M/2})\cup(z_n^-+B_{M/2}).
		\]
		Set $x_n^\pm=t_nz_n^\pm$. Since $t_n\to1$,
		$|x_n^+-x_n^-|\to\infty$ and
		\[
		t_nG_n=(x_n^++B_{m_n/2})\cup(x_n^-+B_{m_n/2}).
		\]
		Moreover,
		\[
		\frac{|t_nF_n\triangle t_nG_n|}{m_n}
		=\frac{|F_n\triangle G_n|}{M}\longrightarrow0.
		\]
		This proves \eqref{eq:two-ball-convergence}.
		
		Since the second term in \eqref{eq:quotient-decomposition}
		is nonnegative, the quotient is bounded below by the first
		term. By $\mathcal A(F_n)\to1$ and the perimeter limit,
		the latter converges to
		\[
		\frac{2P(B_{M/2})-P(B_M)}{M^{2/3}}
		=a(2^{1/3}-1).
		\]
		Hence \eqref{eq:sequence-liminf} follows.
	\end{proof}
	
	\begin{proof}[Proof of Theorem~\ref{thm:main}]
		The coarse bounds follow from Section~\ref{sec:scaling}.
		To prove \eqref{eq:sharp-limit}, let $(m_n)$ be any sequence
		such that $0<m_n<M$ and $m_n\to M$, and set
		$\delta_n=1-m_n/M$.
		By the definition of $\Lambda(m_n)$, choose
		$E_n$ with $|E_n|=m_n$, $P(E_n)<\infty$, and
		$\mathcal A(E_n)>0$ such that
		\begin{equation}\label{eq:near-infimum-choice}
			\Lambda(m_n)\leq
			\frac{\mathcal E_1(E_n)-f(m_n)}
			{m_n^{2/3}\mathcal A(E_n)^2}
			\leq\Lambda(m_n)+\frac{\delta_n}{n}.
		\end{equation}
		The upper bound in \eqref{eq:coarse-bounds}, together with
		\eqref{eq:near-infimum-choice}, gives
		\eqref{eq:bounded-quotient} with a constant independent of $n$.
		
		Moreover,
		\[
		\frac{\Lambda(m_n)}{\delta_n}
		\geq
		\frac{\mathcal E_1(E_n)-f(m_n)}
		{\delta_n m_n^{2/3}\mathcal A(E_n)^2}
		-\frac1n.
		\]
		
		By Proposition~\ref{prop:bounded-quotients}, we get
		\[
		\liminf_{n\to\infty}\frac{\Lambda(m_n)}{\delta_n}
		\geq a(2^{1/3}-1).
		\]
		The upper bound in \eqref{eq:coarse-bounds} also gives
		\[
		\limsup_{n\to\infty}\frac{\Lambda(m_n)}{\delta_n}
		\leq a(2^{1/3}-1).
		\]
		Since $(m_n)$ is arbitrary, \eqref{eq:sharp-limit} follows.
	\end{proof}
	
	\begin{proof}[Proof of Theorem~\ref{thm:almost-optimal}]
		Set $t_n=(m_n/M)^{1/3}$, so that $F_n=t_n^{-1}E_n$.
		By \eqref{eq:almost-optimal} and Theorem~\ref{thm:main},
		the hypothesis of Proposition~\ref{prop:bounded-quotients}
		holds for all sufficiently large $n$. Hence
		\eqref{eq:two-ball-convergence} and
		\eqref{eq:component-limits} follow.
		
		By \eqref{eq:sharp-limit} and \eqref{eq:almost-optimal},
		the left-hand side of \eqref{eq:quotient-decomposition}
		converges to $a(2^{1/3}-1)$. By
		\eqref{eq:component-limits}, the first term on the
		right-hand side has the same limit. Therefore,
		\[
		\frac{t_n^3}{\delta_n}
		\frac{\mathcal E_1(F_n)-f(M)}
		{M^{2/3}\mathcal A(F_n)^2}\to0.
		\]
		Since $t_n^3\to1$ and $\mathcal A(F_n)\to1$,
		we obtain \eqref{eq:critical-small-excess}.
	\end{proof}
	
	\subsection{An equivalent criterion for almost optimality}
	
	\begin{corollary}\label{cor:almost-optimal-criterion}
		Let $m_n\to M$, $0<m_n<M$, $\delta_n=1-m_n/M$, and
		$F_n=(M/m_n)^{1/3}E_n$, where $|E_n|=m_n$ and
		$P(E_n)<\infty$. Assume $\mathcal A(E_n)>0$.
		Then \eqref{eq:almost-optimal} is equivalent to
		\begin{equation}\label{eq:critical-criterion}
			\liminf_{n\to\infty}\mathcal A(F_n)>0,
			\qquad
			\mathcal E_1(F_n)-f(M)=o(\delta_n).
		\end{equation}
	\end{corollary}
	
	\begin{proof}
		The forward implication follows from
		Theorem~\ref{thm:almost-optimal}.
		Conversely, \eqref{eq:critical-criterion} gives a critical
		minimizing sequence and excludes the single-ball alternative.
		Proposition~\ref{prop:critical-dichotomy}, applied to every subsequence,
		then yields the following limits along the full sequence:
		\[
		\mathcal A(F_n)\to1,\qquad
		P(F_n)\to2P(B_{M/2}).
		\]
		By \eqref{eq:quotient-decomposition} and
		\eqref{eq:critical-criterion},
		\[
		\frac{1}{\delta_n}
		\frac{\mathcal E_1(E_n)-f(m_n)}
		{m_n^{2/3}\mathcal A(E_n)^2}
		\longrightarrow a(2^{1/3}-1).
		\]
		Together with \eqref{eq:sharp-limit}, this gives
		\eqref{eq:almost-optimal}.
	\end{proof}
	
	\subsection{Recovery sequences and the critical endpoint}
	
	\begin{corollary}\label{cor:recovery}
		Let $m_n\to M$ with $0<m_n<M$, and set
		$\delta_n=1-m_n/M$. Assume that
		$L_n>2(R_{m_n}+R_{m_n/2})$ for all sufficiently large $n$.
		For
		$E_n=B_{m_n/2}\cup(L_ne_1+B_{m_n/2})$,
		\eqref{eq:almost-optimal} holds if and only if
		\begin{equation}\label{eq:recovery-separation}
			\delta_nL_n\longrightarrow\infty.
		\end{equation}
	\end{corollary}
	
	\begin{proof}
		Lemma~\ref{lem:two-balls} gives
		\begin{equation}\label{eq:recovery-quotient}
			\frac{\mathcal E_1(E_n)-f(m_n)}
			{m_n^{2/3}\mathcal A(E_n)^2}
			=a(2^{1/3}-1)\delta_n+\frac{m_n^{4/3}}{4L_n}.
		\end{equation}
		By \eqref{eq:sharp-limit} and
		\eqref{eq:recovery-quotient}, the difference from $\Lambda(m_n)$
		is $o(\delta_n)$ if and only if
		\[
		\frac{m_n^{4/3}}{4L_n}=o(\delta_n).
		\]
		Since $m_n\to M>0$, this is equivalent to
		$\delta_nL_n\to\infty$.
	\end{proof}
	
	\begin{corollary}
		Define $\Lambda(M)$ by \eqref{eq:optimal-constant} with $m=M$.
		Then $\Lambda(M)=0$, and the infimum is not attained by any set with
		positive asymmetry. The local estimate
		\eqref{eq:critical-local-stability} nevertheless holds.
	\end{corollary}
	
	\begin{proof}
		Proposition~\ref{prop:critical-input} gives $\Lambda(M)\geq0$.
		The competitors in Lemma~\ref{lem:two-balls}, with $m=M$,
		yield $\Lambda(M)\leq0$ as $L\to\infty$. Hence
		$\Lambda(M)=0$.
		The equality case in Proposition~\ref{prop:critical-input}
		excludes attainment by any set with positive asymmetry.
		The final assertion follows from
		Lemma~\ref{lem:local-stability}.
	\end{proof}
	
	\begin{remark}
		Theorems~\ref{thm:main} and \ref{thm:almost-optimal}
		establish the first-order asymptotics of $\Lambda(m)$
		and the qualitative convergence of almost optimal sets.
		They do not imply the exact identity
		$\Lambda(m)=a(2^{1/3}-1)(1-m/M)$ for fixed $m<M$,
		nor a power-law rate in \eqref{eq:two-ball-convergence}.
		The separation condition \eqref{eq:recovery-separation}
		applies only to the explicit two-ball family. The qualitative
		profile decomposition yields no corresponding rate for
		general almost optimal sequences.
		The truncation bounds $4/n$ and $8/n^2$ in
		\eqref{eq:profile-approximation} are used only
		in the critical $o(1)$ compactness argument.
		
		For the family in Corollary~\ref{cor:recovery},
		choose $L_n=\delta_n^{-1}$ for all sufficiently large $n$
		and set $F_n=(M/m_n)^{1/3}E_n$.
		Then $\mathcal A(E_n)=1$, and the symmetric difference
		in \eqref{eq:two-ball-convergence} vanishes.
		By \eqref{eq:two-ball-exact}, applied at mass $M$,
		\[
		\begin{aligned}
			\frac{\mathcal E_1(F_n)-f(M)}{\delta_n}
			=\,\frac{M^2}{4\delta_n(M/m_n)^{1/3}L_n}
			=\,\frac{M^2}{4}\left(\frac{m_n}{M}\right)^{1/3}
			\longrightarrow\frac{M^2}{4}>0.
		\end{aligned}
		\]
		Thus the energy condition in \eqref{eq:critical-criterion}
		is not satisfied, distinguishing this family from
		almost optimal sequences.
	\end{remark}
	
	\section*{Acknowledgments}
     The research of Chengxiang Zhang is supported by the National Natural Science Foundation of China (Grant No. 12371107). The authors acknowledge the use of ChatGPT (OpenAI) to assist with
	selected calculations. They have independently verified all
	mathematical arguments and calculations and take full responsibility
	for the content of the paper.

\end{document}